\documentclass[a4paper, 12pt]{article}

\usepackage[koi8-r,cp1251]{inputenc}
\usepackage[english,russian]{babel}
\usepackage{amsfonts,amssymb,amsmath, hyperref}
\usepackage[final]{epsfig}
\usepackage{graphicx}

\newtheorem{theorem}{Теорема}
\newtheorem{lemma}{Лемма}
\newtheorem{proposition}{Предложение}

\newtheorem{definition}{Определение}
\newtheorem{corollary}{Следствие}

\begin{document}

УДК 517.986

\begin{center}
\Large{\textbf{О границе Шилова и спектре Гельфанда алгебр
обобщенных аналитических функций}}
\end{center}

\begin{center}
А.Р.\,Миротин
\end{center}

\begin{center}
amirotin@yandex.ru
\end{center}




В работе изучается алгебра обобщенных аналитических
функций, определенных на полугруппе полухарактеров $\widehat {S}$
дискретной абелевой полугруппы $S$ с сокращениями и единицей,
более широкая, чем алгебра Аренса-Зингера. Показано, что строгая
граница и граница Шилова этой алгебры являются объединениями
максимальных подгрупп полугруппы $\widehat {S}$. Если $S$ не
содержит нетривиальных простых идеалов, обе границы совпадают с
группой характеров полугруппы $S$. В последнем случае вычислен и
спектр Гельфанда рассматриваемой алгебры.

\

Ключевые слова.
Граница Шилова, спектр Гельфанда, равномерная алгебра,
обобщенная аналитическая функция.

\

Abstract.
Let $S$ be a discrete abelian semigroup with unit
and concellations and $\widehat {S}$ the
semigroup of semicharacters of  $S$. We shaw that the strong boundary and the Shilov
boundary of the algebra of generalized analytic functions defined on $\widehat {S}$  are
 unions of some maximal subgroups of
$\widehat {S}$.  We shaw also that the both boundaries coincide with the character group of $S$ if $S$ does not contain nontrivial simple ideals. In
the last case the Gelfand spectrum of the algebra under consideration
is calculated.

\
Keywords. Shilov boundary, Gelfand spectrum, uniform
algebra, generalized analytic function.

\



 \section{Введение}
    Первые работы по теории обобщенных аналитических
 функций на пространствах полухарактеров принадлежат Р.~Аренсу
и И.М.~Зингеру \cite{AS}, \cite{A}. В настоящее время это довольно
развитый раздел функционального анализа, расположенный на стыке
теории равномерных алгебр, абстрактного гармонического анализа и
теории аналитических функций одного и нескольких комплексных
переменных. Современное состояние предмета отражено в монографиях
\cite{G} -- \cite{GT}.

С другой стороны, имеются задачи, для решения которых понятие
обобщенной аналитичности в смысле Аренса-Зингера оказывается
слишком ограничительным. Так, в работе \cite{M1} была решена
задача описания образа пространства квадратично суммируемых
функций на подполугруппе $S$ топологической абелевой группы
относительно преобразования Лапласа на этой полугруппе. В случае
полуоси ответ на этот вопрос дает, как известно, классическая
теорема Пэли-Винера. Согласно этой теореме образом в этом случае
служит пространство Харди $H^2$ в полуплоскости, которое состоит
из функций, аналитических в этой полуплоскости и удовлетворяющих
некоторым дополнительным ограничениям. В случае общих полугрупп,
рассмотренном в \cite{M1}, ответом служит обобщенное пространство
Харди $H^2(\widehat S)$, состоящее из  функций, определенных на
пространстве $\widehat S$ полухарактеров  полугруппы $S$ и
аналитических в смысле, более общем, чем в \cite{AS}. В то же
время, равномерная алгебра $A(\widehat S)$, состоящая из функций,
непрерывных на пространстве полухарактеров дискретной полугруппы и
аналитических "внутри"\  этого пространства в этом более общем,
чем у Аренса и Зингера, смысле, оставалась неисследованной. Такие
исследования полезны и для самой теории Аренса-Зингера. Например,
представляет интерес вопрос об условиях, при которых алгебра
Аренса-Зингера и алгебра $A(\widehat S)$ совпадают, поскольку его
решение  дает новое удобное определение аналитичности по
Аренсу-Зингеру. Другая задача теории алгебр Аренса-Зингера, при
решении которой возникают алгебры типа $A(\widehat S)$,
--- это проблема описания интерполяционных подмножеств (см. \cite{VUZ07}).
Сужения  алгебр типа $A(\widehat S)$ на группу характеров
полугруппы $S$ дают новые примеры сдвигово-инвариантных алгебр в
смысле (\cite{GT}, с. 120).  Добавим, что для многих полугрупп
соответствующее пространство $\widehat S$ является (вещественным
или комплексным) многообразием, что приводит к результатам об
алгебрах аналитических функций на этих многообразиях.

 Следует заметить, что теория алгебр типа $A(\widehat S)$ оказалась
в некоторых отношениях более трудной, чем теория алгебр
Аренса-Зингера, и в ней возникает ряд интересных задач. Например,
проблема описания спектра Гельфанда таких алгебр в общем случае
остается нерешенной.

Основные результаты данной работы описывают структуру строгой
границы и границы Шилова равномерной алгебры функций,
аналитических в смысле \cite{M1}. Если полугруппа не содержит
нетривиальных простых идеалов, обе границы совпадают с группой
характеров этой полугруппы. В последнем случае получен ответ на
 вопрос о  спектре Гельфанда
рассматриваемой алгебры. При других предположениях эти вопросы
рассматривались в \cite{VUZ07}. Отметим также, что нижеследующие
теоремы 1 и 2 усиливают результаты, анонсированные в \cite{MR}.
Статья продолжает цикл работ автора, посвященных гармоническому анализу на полугруппах,
см. \cite{Mir2}, \cite{Hilb}, \cite{SbMath}.

Всюду ниже $S$ --- дискретная абелева полугруппа с сокращениями и
единицей $e$, записываемая мультипликативно, $G=S^{-1}S$ ---
группа частных для $S$.
    {\it Полухарактером} полугруппы $S$ называется
гомоморфизм $\psi$ полугруппы $S$ в мультипликативную полугруппу
$\overline {\Bbb{D}}=\{ z\in~\Bbb{C}: |z|\leq 1\}$, не являющийся
тождественным нулем.  {\it Характерами}  называются полухарактеры,
равные по модулю единице.

    Множество всех полухарактеров полугруппы $S$ далее обозначается
    через $\widehat {S}$,    а его подмножество,
состоящее из неотрицательных полухарактеров, --- через $\widehat
{S}_+$. Символом $E(\widehat {S})$ будем обозначать множество всех
{\it идемпотентов} полугруппы $\widehat {S}$, т. е.
полухарактеров, принимающих лишь значения 0 и 1. Множества
$\widehat {S}, \widehat {S}_+$ и  $E(\widehat {S})$, наделенные
топологией поточечной сходимости на $S$, являются компактными
топологическими полугруппами по умножению с единицей $1$ (они
замкнуты в $\overline {\Bbb{D}}^S$). Компактную группу всех
характеров полугруппы $S$ будем обозначать $X$.

 Согласно   (\cite{AS}, теорема 3.1) каждый полухарактер
$\psi\in\widehat {S}$ может быть представлен в виде
$\psi=\rho\chi$, где $\chi\in X$, а полухарактер $\rho=|\psi|\in
\widehat {S}_+$ определяется по $\psi$ однозначно ({\it полярное
разложение}).
 С  полухарактером $\rho\in \widehat {S}_+$ связаны
подполугруппы $S(\rho):=\{s\in S: \rho(s)> 0\}$ и
$S^{\rho}:=\{s\in S: \rho(s)=1\}$ полугруппы $S$, дополнения
которых, если они не пусты, являются идеалами полугруппы $S$.
Идеалы, дополнения которых есть полугруппа, называются {\it
простыми}. Простые идеалы, отличные от $S\setminus\{e\}$, будем
называть {\it нетривиальными}. Через $\theta_S$ (или просто
$\theta$, если это не будет приводить к недоразумениям) обозначим
индикатор подмножества $\{e\}\subset S$. Это полухарактер
полугруппы $S$, если и только если $S^{-1}\cap S=\{e\}$. Наконец
отметим, что степень $\rho^0$ по определению есть индикатор
$S(\rho)$, и что $\rho^z\in \widehat {S}\setminus X$ при $\rho\in
\widehat {S}_+, \quad \rho\ne 1,\quad z\in \Pi$, где $\Pi:=\{ {\rm
Re} z>0\}$ (\cite{AS}, \S 7). Относительно понятий и обозначений
из теории коммутативных банаховых алгебр, не определенных в данной
работе, см. \cite{G}.

 \section{Границы}
Пусть $A$ есть коммутативная банахова алгебра со спектром
Гельфанда  $M_A$. {\it Точечной производной} в точке $\phi\in M_A$
называется такой линейный функционал $L$ на $A$, для которого
$$
L(fg)=\phi(f)L(g)+\phi(g)L(f),\quad f,g\in A.
$$
\noindent В частности, если $A$ есть равномерная алгебра на
компакте $Y$, то точечная производная в точке $y\in Y$
удовлетворяет равенству
$$
L(fg)=f(y)L(g)+g(y)L(f),\quad f,g\in A.
$$

Напомним, что {\it строгой границей} равномерной алгебры $A$ на
компакте $Y$ называется множество  всех ее $p$-точек, т. е. точек
из $Y$, являющихся пересечениями множеств пика.

Для доказательства основных результатов нам понадобятся три леммы.
 Следующее утверждение сформулировано в (\cite{G}, c. 92). Ввиду
 отсутствия ссылок на работу, где эта лемма доказана, приводим
 здесь ее доказательство.

\begin{lemma}
Пусть $A$ есть равномерная алгебра на компакте $Y$. Если точка
$p\in Y$ принадлежит строгой границе алгебры $A$, то любая
точечная производная в точке $p$ тривиальна.
\end{lemma}

Доказательство. Рассмотрим максимальный идеал $I_p=\{f\in
A:f(p)=0\}$ алгебры $A$ и покажем, что он имеет ограниченную
аппроксимативную единицу. Ввиду теоремы Альтмана (см., например,
\cite{BD}, c. 58) достаточно проверить, что для любых $f\in I_p$ и
$\epsilon>0$ найдется такой $u\in I_p$, $\|u\|\leq 2$, что
$\|f-fu\|<\epsilon$. Для проверки выберем  окрестность $U$ точки
$p$, для которой $|f(x)|<\epsilon$ при $x\in U$ . В силу
(\cite{G}, лемма II.12.2)  найдется такое множество пика  $F$, что
$p\in F\subset U$. Если функция $h$ образует пик на $F$, то для
натурального $n$, при котором $|h(y)|^n<\epsilon/\|f\|$ для $y\in
Y\setminus U$, функция $u=1-h^n$ удовлетворяет поставленным
требованиям. Теперь в силу факторизационной теоремы Коэна (см.,
например, \cite{HR}, следствие (32.26)) для любого $f\in I_p$
найдутся такие $f_1, f_2\in I_p$, что $f=f_1f_2$. Если $L$ есть
точечная производная в точке $p$, то
$L(f)=f_1(p)L(f_2)+f_2(p)L(f_1)=0$. Поскольку коразмерность $ I_p$
равна 1, и $L(1)=0$, то $L=0.\Box$

Следующее понятие в неявной форме впервые появилось в \cite{M1}, а
в явной -- в \cite{VUZ07}.

\begin{definition}
  Комплекснозначную
функцию $F$ на $\widehat {S} \setminus X$ будем называть {\it
обобщенной аналитической}, если для любых полухарактеров $\rho,
\psi$ из $\widehat {S}\setminus X$, $\rho \geq 0$ отображение $z
\mapsto F(\rho^{z} \psi)$
  аналитично в открытой правой полуплоскости $\Pi$ и непрерывно в +0.
\end{definition}

Равномерную алгебру функций, непрерывных на $\widehat {S}$ и
обобщенных аналитических на $\widehat {S}\setminus X$, обозначим
$A(\widehat {S})$.  Границу Шилова алгебры $A(\widehat {S})$ будем
обозначать через $\partial_{A(\widehat {S})}$, а строгую границу
 -- через $\Gamma$.

Легко проверить, что функции вида  $\widehat a(\xi):=\xi(a)\ (a\in
S,\ \xi\in\widehat {S})$ принадлежат $A(\widehat {S})$.
Равномерная алгебра  на $\widehat {S}$, порожденная этими
функциями, называется алгеброй Аренса-Зингера и обозначается
$A_0(\widehat {S})$. Ясно, что $A_0(\widehat {S})\subseteq
A(\widehat {S})$,  причем строгое включение возможно
(\cite{VUZ07}, пример 1).

\begin{lemma}
Если полухарактер $\psi$  полугруппы $S$ принадлежит строгой
границе $\Gamma$ алгебры $A(\widehat {S})$, то $|\psi|\in
E(\widehat {S})$.
\end{lemma}

Доказательство. Пусть $\psi=\rho\chi$ -- полярное разложение
полухарактера $\psi$. Допустим противное, т. е. $\rho\notin
E(\widehat {S})$. Если  положим $\kappa(z)=\rho^{z} \psi\
(z\in\Pi)$, то функционал, определенный на функциях $F\in
A(\widehat {S})$ равенством $L(F)=(F\circ\kappa)'(1),$
 есть, как легко проверить, точечная производная алгебры
$A(\widehat {S})$ в точке $\psi$. Кроме того, если элемент $a\in
S$ таков, что  $\rho(a)\ne 0,1$, то $L(\widehat
a)=\psi(a)\log\rho(a)\ne 0$. Применение леммы 1 приводит  к
противоречию.$\Box$

Обозначим через $E_0(\widehat {S})$ множество тех полухарактеров
$\rho_0\in E(\widehat {S})$, для каждого из которых найдется такой
полухарактер $\rho_1\in \widehat {S}_+$, что $\rho_1(s)=1$  при
$s\in S(\rho_0)$  и  $0<\rho_1(s)<1$  при $s\notin S(\rho_0)$.

Нижеследующее предложение показывает, что утверждение, обратное
лемме 2, вообще говоря, неверно.

 \begin{proposition}
$\theta\notin E_0(\widehat {S})\cap\Gamma$.
\end{proposition}

Доказательство. Пусть $\theta\in E_0(\widehat {S})$. В силу
  (\cite{VUZ07}, теорема 2) для любой функции $F\in A(\widehat {S})$
справедливо равенство $\int\limits_X F(\chi)d\chi=F(\theta),$
 где $d\chi$  есть нормированная  мера Хаара (компактной)
группы $X$. Таким образом, комплексный гомоморфизм $F\mapsto
F(\theta)$ алгебры $A(\widehat {S})$ имеет представляющую меру
$\sigma\ne\delta_\theta$, что невозможно для точек строгой границы
(см., например, \cite{G}, теорема II.11.3 и замечание  с. 86;
\cite{Tay}, утверждение 10.3.1).$\Box$

Ниже (см. теорему 1)  обобщим предложение 1 для случая, когда
полугруппа $S$ счетна.

В связи со следующей леммой отметим, что множество $E(T)$
идемпотентов абелевой компактной полугруппы $T$ есть компактная
полугруппа, и для каждого $e\in E(T)$ существует наибольшая
подгруппа $H(e)$, содержащая $e$.  Эта группа компактна и
называется {\it максимальной подгруппой полугруппы $T$, содержащей
$e$}; кроме того,  множество  $H=\coprod_{e\in E(T)} H(e)$ также
компактно в $T$ (см., например, \cite{CHC}, c. 17 - 18).

\begin{lemma}
 Для любого полухарактера $\rho_0\in
E(\widehat {S})$ справедливо равенство  $H(\rho_0)=\rho_0X$.
\end{lemma}

Доказательство. Поскольку $\rho_0X$  является группой, содержащей
$\rho_0$, то $\rho_0X\subseteq H(\rho_0)$. Пусть теперь $\psi\in
H(\rho_0),\ |\psi|=\rho$. Так как $\rho_0\psi=\psi$, то, переходя
к модулям, получаем $\rho_0\rho=\rho$. Следовательно,
$S(\rho)\subseteq S(\rho_0)$. С другой стороны, найдется такой
$\psi'\in H(\rho_0)$, что $\psi\psi'=\rho_0$. Если
$\rho'=|\psi'|$, то  $\rho\rho'=\rho_0$, а потому $\rho=\rho_0.
\Box$

\begin{theorem}
  Строгая граница $\Gamma$ алгебры
 $A(\widehat {S})$
 есть объединение максимальных подгрупп полугруппы $\widehat {S}$,
 т. е.
$$
\Gamma=\coprod\limits_{\rho\in F}\rho X,
$$
\noindent где $F$  содержится в  $E(\widehat {S})$. Если $S$
счетна, то $1\in F$, а $E_0(\widehat {S})\cap F=\emptyset$.
\end{theorem}

Доказательство. По лемме 2 любой полухарактер $\psi\in \Gamma$
имеет полярное разложение $\psi=\rho\chi$, в котором $\rho\in
E(\widehat {S})$. С учетом леммы 3 это означает, что $\psi$
принадлежит максимальной подгруппе $\rho X$ полугруппы $\widehat
{S}$. Заметим, что каждому характеру $\xi\in X$ соответствует
автоморфизм алгебры $A(\widehat {S})$, переводящий $f\in
A(\widehat {S})$ в $f_\xi:\psi\mapsto f(\xi\psi)$. Отсюда следует,
что $\xi\Gamma=\Gamma$, а потому $\psi X=\rho X\subseteq \Gamma$.
Если $F$ есть образ $\Gamma$ при отображении $\psi\mapsto |\psi|$,
то $ \Gamma=\coprod_{\rho\in F}\rho X$.

Далее будем  предполагать, что $S$ счетна, $S=\{s_n:
n=1,2,\ldots\}$. Зафиксируем  строго положительную функцию $v$ на
$S$, удовлетворяющую условию $\sum_{n=1}^\infty v(s_n)=1$, а также
точку $\chi\in X$ и рассмотрим функцию
$$
f(\psi):=\sum\limits_{n=1}^\infty\psi(s_n)\overline{\chi(s_n)}v(s_n)
=\langle\psi,\chi\rangle,
$$
\noindent где угловые скобки обозначают скалярное произведение в
весовом пространстве $l_2(S,v)$. Ясно, что $f\in A(\widehat {S})$
как сумма равномерно по $\psi$ сходящегося ряда аналитических
функций (мажорантным рядом является $\sum_{n=1}^\infty v(s_n)=1$).
Неравенство Коши-Буняковского $|f(\psi)|\leq\|\psi\|\|\chi\|\leq
1$ показывает, что $|f|$ достигает своего максимума, равного $1$,
в точке  $\chi$ и только в этой точке (если оно при некотором
$\psi\in \widehat {S}$ превращается в равенство, то
$\psi(s)=c\chi(s)$, и при $s=e$ получаем $c=1$). Таким образом,
$X$  содержится в любой границе алгебры $A(\widehat {S})$, а тогда
$1\in F$.

Заметим, наконец, что по критерию А. С. Мищенко (см., например,
\cite{K}, c. 391) компактное пространство $\widehat {S}$
метризуемо, так как каждая его точка $\psi$ обладает счетной базой
окрестностей вида
$$
U(\psi; t_1,\ldots,t_n; \epsilon):=\{\zeta\in\widehat
{S}:|\zeta(t_i)-\psi(t_i)|<\epsilon \mbox { при } i=1,\ldots,n \},
$$
\noindent где $t_i\in S, n\in \mathbb{N}, \epsilon\in
\mathbb{Q}_+$.  В этой ситуации, как известно, $\Gamma$ есть
наименьшая из границ алгебры $A(\widehat {S})$. Пусть $f\in
A(\widehat {S}), \rho_0\in E_0(\widehat {S})$ и  $\rho_1\in
\widehat {S}_+$ таков, что $\rho_1(s)=1$  при $s\in S(\rho_0)$  и
$0<\rho_1(s)<1$  при $s\notin S(\rho_0)$. Так как функция
$\phi(z)=f(\rho_1^z)$ непрерывна и ограничена в замкнутой правой
полуплоскости $\overline{\Pi}$ и аналитична в $\Pi$, то по
принципу Фрагмена-Линделефа $\sup_{z\in
\overline{\Pi}}|\phi(z)|=\sup_{r\in \mathbb{R}}|\phi(ir)|$. Отсюда
следует, что при всех натуральных $n$ справедливо неравенство
$$
|f(\rho_1^n)|\leq\sup\limits_{r\in \mathbb{R}}|f(\rho_1^{ir})|.
$$
\noindent Поскольку $\rho_1^{ir}\in X$, то, устремляя $n$ к
$\infty$, получаем $|f(\rho_0)|\leq\max\limits_X |f|$. Но
$X\subset \Gamma$, и $\Gamma$ минимальна. Следовательно,
$\rho_0\notin \Gamma$. $\Box$

Отметим, что в условиях предложения  полухарактер $\theta\notin
F$.

Выясним теперь структуру границы Шилова алгебры $A(\widehat {S})$
(ниже квадратные скобки обозначают замыкание в $\widehat {S}$).

\begin{theorem}
 Для границы Шилова  алгебры $A(\widehat {S})$
 справедливо равенство
$$
\partial_{A(\widehat
{S})}=\coprod\limits_{\rho\in K}\rho X,
$$
\noindent где $K=[F]$  есть компактное подмножество  $E(\widehat
{S})$, содержащее $1$.
\end{theorem}

Доказательство. Заметим сначала, что для любого замкнутого
$C\subset E(\widehat {S})$ множество $M:=\coprod_{\rho\in C}\rho
X$ замкнуто в $\widehat {S}$. Действительно, если направленность
$\psi_n$ из $M$ сходится к полухарактеру $\psi\in  \widehat {S}$,
то $|\psi_n|$  сходится к  $|\psi|$. А так как $|\psi_n|\in C$, то
и $|\psi|\in C$. В силу полярного разложения $\psi\in|\psi|\cdot
X\subset M$.

Покажем теперь, что для любого множества $D\subset E(\widehat
{S})$ справедливо равенство
$$
\left[\coprod\limits_{\rho\in D}\rho
X\right]=\coprod\limits_{\rho\in [D]}\rho X.
$$

  Действительно, левая часть содержится в правой по причине
замкнутости последней. Для доказательства обратного включения
возьмем полухарактер $\psi\in \coprod_{\rho\in [D]}\rho X$. Он
имеет полярное разложение $\psi=\rho\chi$, где $\rho\in [D],\
\chi\in X$. Тогда $\rho=\lim_n\rho_n$  для некоторой
направленности $\rho_n\in D$. Поэтому направленность
$\psi_n:=\rho_n\chi\to \psi$, причем  $\psi_n\in
\coprod\limits_{\rho\in D}\rho X$. Значит, $\psi\in
\left[\coprod\limits_{\rho\in D}\rho X\right]$.

Теперь, поскольку  строгая граница $\Gamma$ плотна в
$\partial_{A(\widehat {S})}$  (см., например, \cite{Tay}, c. 90),
с учетом теоремы 1 имеем
$$
\partial_{A(\widehat {S})}=[\Gamma]=\coprod\limits_{\rho\in [F]}\rho
X.
$$

Компактность  $[F]$  следует из компактности $E(\widehat {S})$.
Наконец, заметим, что алгебра $A(\widehat {S})$ содержит в
качестве замкнутой подалгебры алгебру Аренса-Зингера $A_0(\widehat
{S})$, границей Шилова которой служит $X$ \cite{AS}, а потому
$X\subset\partial_{A(\widehat {S})}$.$\Box$

\begin{corollary}
 Если $E(\widehat {S})$ содержит лишь конечное
 множество $p$-точек,
 то $\Gamma =\partial_{A(\widehat {S})}$. Кроме того,
 $\theta\notin E_0(\widehat {S})\cap \partial_{A(\widehat {S})}$.
\end{corollary}

Доказательство. По условию множество $F=\Gamma\cap E(\widehat
{S})$  конечно. Поэтому первое утверждение следует из равенства
$K=[F]$, а второе --- из предложения 1. $\Box$

\begin{corollary}
 Если полугруппа  $S$ не содержит
нетривиальных простых идеалов, то
$$
\Gamma =\partial_{A(\widehat {S})}=X.
$$
\end{corollary}

Доказательство. В силу теоремы 2 имеем $X\subset
\partial_{A(\widehat {S})}$. Поскольку дополнения к $S(\rho)$ и
$S^{\rho}$ в случае их непустоты  являются простыми идеалами, а
множество $S\setminus\{e\}$ есть по условию единственный простой
идеал полугруппы $S$, то все неотрицательные полухарактеры $\rho$,
отличные от $\theta$, удовлетворяют условию $0<\rho(s)<1$ при всех
$s\in S, s\ne e$, а  идемпотентами полугруппы $\widehat {S}$
являются только $1$ и $\theta$
 ($S\setminus(S^{-1}\cap S)$ есть простой идеал, и значит,
 $S^{-1}\cap S=\{e\}$).
  Поэтому в  силу следствия
1 и теоремы 1 $\Gamma =\partial_{A(\widehat {S})}\subseteq
X\cup\{\theta\}$.   Из (\cite{Tay},  теорема 4.5.2) следует
$\widehat {S}_+\ne\{1,\theta\}$ (см. также \cite{M3}, теорема
2.3.9), а значит выполнены все условия следствия 1 и $\theta\in
E_0(\widehat S)$. Следовательно, $\theta\notin\partial_{A(\widehat
{S})}.$ $ \Box$

\section{ Спектр Гельфанда}
Следствие 2 позволяет вычислить спектр Гельфанда алгебры
$A(\widehat {S})$ в случае, когда  $S$ не содержит нетривиальных
простых идеалов. Ниже через $A(X)$ обозначаем алгебру,
образованную сужениями на $X$ функций из $A(\widehat {S})$, а
$S_1$ есть подполугруппа группы $G$,
 порожденная носителями преобразований Фурье
функций из  $A(X)$.

\begin{theorem}
 Пусть полугруппа  $S$ не содержит нетривиальных простых идеалов.
Тогда

{\rm 1)} $S_1^{-1}\cap S_1=\{e\}$;

{\rm 2)} $\widehat {S}$ есть подполугруппа полугруппы $\widehat
{S_1}$ в том смысле, что $S\subseteq S_1$ и отображение сужения на
$S$ есть топологический изоморфизм некоторой подполугруппы
полугруппы $\widehat {S_1}$ и полугруппы $\widehat {S}$;

{\rm 3)} отображение сужения на $\widehat {S}$ есть изометрический
изоморфизм алгебры Аренса-Зингера $A_0(\widehat {S_1})$ и алгебры
$A(\widehat {S})$;

{\rm 4)} спектр Гельфанда алгебры $A(\widehat {S})$ можно
отождествить с $\widehat {S_1}$.
\end{theorem}

Доказательство. Сначала докажем 1, 2).  Для любого $a\in S$
обозначим через $\widetilde{a}$ функцию вычисления $\psi\mapsto
\psi(a)\ (\psi\in \widehat {S})$. Тогда $\widetilde{a}\in
A(\widehat{S})$. Заметим, прежде всего, что   носителем
преобразования Фурье сужения функции $\widetilde{a}$ на $X$ служит
множество $\{a\}$. Поэтому $S\subseteq S_1$. Также легко
проверить, что группы характеров полугрупп $S$ и $S_1$ естественно
изоморфны группе характеров $X$ группы $G$ (характеры этих
полугрупп однозначно продолжаются на $G$), поэтому далее  будем их
отождествлять.

Поскольку алгебра $A(X)$, очевидно,  инвариантна относительно
 сдвигов, она в
силу (\cite{GT}, предложение 4.1.8) равна алгебре $A_{S_1}$,
состоящей из всех непрерывных на $X$ функций, преобразование Фурье
которых сосредоточено на $S_1$. Последняя алгебра, в свою очередь,
посредством  отображения сужения  $i_0:F\mapsto F|X$ изометрически
изоморфна алгебре $A_0(\widehat {S_1})$ по теореме 2.6 из
\cite{AS}. С другой стороны, равенство $\partial_{A(\widehat
{S})}=X$ влечет то, что отображения сужения  $i:\Phi\mapsto
\Phi|X$ есть изометрический изоморфизм алгебр $A(\widehat {S})$ и
$A(X)$. Таким образом, возникает изометрический изоморфизм
$j:A_0(\widehat {S_1})\to A(\widehat {S})$ такой, что следующая
диаграмма коммутативна:

\begin{displaymath}
\begin{array}{ccc}
A_{S_1} & = & A(X) \\
\uparrow\lefteqn{i_0}&&\uparrow\lefteqn{i}\\
A_0(\widehat {S_1})&\stackrel{j}{\longrightarrow}& A(\widehat{S})
\end{array}
\end{displaymath}

Покажем, что каждый полухарактер $\psi\in\widehat{S}$ продолжается
до полухарактера $\widetilde{\psi}\in \widehat {S_1}$.
Действительно, для любого $x\in S_1$ имеем $\widehat {x}\in
A_0(\widehat {S_1})$. Положим $\widetilde{\psi}(x):=j(\widehat
{x})(\psi)$. Тогда при всех $x,y\in S_1$
$$
\widetilde{\psi}(xy)=j(\widehat {xy})(\psi)=j(\widehat
{x})(\psi)j(\widehat
{y})(\psi)=\widetilde{\psi}(x)\widetilde{\psi}(y).
$$
Кроме того, изометричность $j$ дает
$|\widetilde{\psi}(x)|=|j(\widehat {x})(\psi)|\leq \|j(\widehat
{x})\|=\|\widehat {x}\|=1$.  Следовательно, $\widetilde{\psi}\in
\widehat {S_1}$.

Покажем, что $\widetilde{\psi}|S=\psi$. Заметим, что
$j(\widehat{a})=i^{-1}(i_0(\widehat{a}))=i^{-1}(\widehat{a}|X)
=\widetilde{a}$ ($a\in S$), поскольку функции из $A(\widehat{S})$
и $A_0(\widehat {S_1})$  однозначно определяются своими значениями
на $X$, а значения $\widehat{a}$ и $\widetilde{a}$ на $X$
совпадают. Поэтому для любого  $\psi\in \widehat {S_1}$ имеем
$\widetilde{\psi}(a)=j(\widehat{a})(\psi)=\widetilde{a}(\psi)=\psi(a)$,
что и утверждалось.

Как показано в доказательстве следствия 2, если полухарактер
$\psi\ne\theta_S$, то $\psi(s)\in (0,1)$  при $s\ne e$.
Следовательно, если $x=a^{-1}b\in S_1\ (a,b\in S)$, то необходимо
$\widetilde{\psi}(x)=\psi(b)/\psi(a)$. Теперь ясно, что (очевидно,
непрерывное и инъективное) отображение $h:\psi\mapsto
\widetilde{\psi}$  из $\widehat{S}$  в $\widehat{S_1}$ есть
гомоморфизм полугруппы  $\widehat{S}\setminus\{\theta_S\}$  в
$\widehat {S_1}$. В доказательстве следствия 2 показано также, что
$S^{-1}\cap S=\{e\}$, а потому в силу (\cite{V07}, лемма 3)
полугруппа $\widehat{S}$  связна. Следовательно, полугруппа
$\widehat{S}\setminus\{\theta_S\}$ плотна в $\widehat{S}$, и
значит $h$ однозначно продолжается по непрерывности до
гомоморфизма $h:\widehat{S}\to \widehat{S_1}$. Так как $\widehat
{S}$ компактна, то $h$ есть гомеоморфизм $\widehat {S}$  на свой
образ. Следовательно, отображение $h^{-1}$ (являющееся
отображением сужения на $S$) есть топологический изоморфизм
полугрупп $h(\widehat{S})$ и $\widehat{S}$. Это доказывает 2).

Для завершения доказательства утверждения 1) достаточно установить
равенство  $h(\theta_S)=\theta_{S_1}$, так как из него следует,
что $\theta_{S_1}$ есть полухарактер полугруппы $S_1$. Но
$h(\theta_S)(h(\psi)-1)=0$, поскольку $h$ есть гомоморфизм,  т. е.
$\widetilde{\theta_S}(\widetilde{\psi}-1)=0$   для любого
$\psi\in\widehat{S}$. Пусть $x\in S_1\setminus\{e\}$, и
$\psi=\chi|S$, где $\chi\in X,\ \chi(x)\ne 1$. Тогда равенство
$\widetilde{\theta_S}(x)(\widetilde{\psi}(x)-1)=0$ влечет
$\widetilde{\theta_S}(x)=0$, что и требовалось доказать.

3) Утверждение следует из 2) и того факта, что для любой функции
$F\in A_0(\widehat {S_1})$ функции $j(F)$ и $F|\widehat {S}$
совпадают на  $X$, а потому и всюду.

 4) Это сразу вытекает из утверждения 3) и теоремы 4.1
[1], согласно которой спектр Гельфанда алгебры $A_0(\widehat
{S_1})$ есть $\widehat {S_1}$. $\Box$

{\bf Пример}. Выберем иррациональные числа $\alpha, \beta$, такие,
что $\beta<\alpha<0$, и рассмотрим следующую подполугруппу группы
$\mathbb{Z}^2$: $S_{\alpha\beta}=\{(m,n)\in \mathbb{Z}^2|\alpha
m\leq n, \beta m\leq n\}.$  Полугруппа $S_{\alpha\beta}$ не
содержит нетривиальных простых идеалов, и ее полухарактеры имеют в
точности вид $\psi(m,n)=z_1^m z_2^n$, где $z_1,z_2\in
\overline{\mathbb{D}},\ |z_1|^{\alpha}\leq |z_2|\leq
|z_1|^{\beta}\ (0^0:=1)$. Таким образом, полугруппа
$\widehat{S_{\alpha\beta}}$ естественным образом отождествляется с
замкнутой областью в $\mathbb{C}^2$
$$
\Delta_{\alpha\beta}=\{(z_1,z_2)\in \overline{\mathbb{D}}^2|
|z_1|^{\alpha}\leq |z_2|\leq |z_1|^{\beta}\}.
$$
Из вида полухарактеров следует, что $S_{\alpha\beta}$ является
конусом в смысле \cite{M1}, а потому в силу (\cite{VUZ07}, теорема
3) алгебра $A(\widehat{S_{\alpha\beta}})$ совпадает с алгеброй
Аренса-Зингера $A_0(\widehat{S_{\alpha\beta}})$, отвечающей
полугруппе $S_{\alpha\beta}$.

Предположим теперь, что $\alpha-\beta>1$. Тогда множество
 $S=S_{\alpha\beta}\setminus\{(0,1)\}$ является полугруппой, которая тоже не
содержит нетривиальных простых идеалов. Так как ее полухарактеры,
отличные от $\theta$, не принимают нулевых значений, они
продолжаются на полугруппу $S_{\alpha\beta}$, и  имеют тот же вид,
что и у этой полугруппы. Алгебра $A(\widehat{S})$ состоит из
функций, непрерывных на $\Delta_{\alpha\beta}$ и аналитических
внутри этого множества, и не совпадает с $A_0(\widehat{S})$. В
этом примере $S_1=S_{\alpha\beta}$. Следовательно, строгая граница
и граница Шилова алгебры $A(\widehat{S})$ отождествляются с
двумерным тором $\mathbb{T}^2$, а ее спектр Гельфанда --- с
пространством $\Delta_{\alpha\beta}$.



\end{document}